\documentclass[12pt,reqno]{amsart}
\usepackage{amsmath,amsfonts,  amssymb,amsthm,amscd}
\makeindex

\usepackage{amssymb}
\usepackage{graphics}
\usepackage[latin1]{inputenc}
\usepackage{epsf}

 \usepackage{ epsfig,epsf}
\usepackage{enumerate}
\usepackage{indentfirst}
\usepackage{euscript}
\usepackage{graphicx}
\usepackage{amsmath}
\usepackage{amsfonts}
\usepackage{amssymb}
\makeindex

 \theoremstyle{plain}
\newtheorem{proposition}{Proposition}

\newtheorem{remark}{Remark}
\theoremstyle{definition}
\newtheorem{definition}{Definition}
\newtheorem{theorem}{Theorem}

 \title[   Bifurcations of Umbilic Points]{  Bifurcations of Umbilic Points
 and Related Principal Cycles}

 \author[C. Guti\'errez, J. Sotomayor  and R. Garcia]{Carlos Guti\'errez, \;Jorge Sotomayor \;and \;Ronaldo Garcia}


 \thanks{The  authors are fellows of  CNPq.
 This work was done under the project PRONEX/CNPq/MCT - grant number 66.2249/1997-6
 Teoria Qualitativa das Equa\c c\~oes Diferenciais Ordin\'arias and was
 partially supported by CNPq Grant  476886/2001-5.}
 
\begin{document}
 \maketitle

 {\small {\sc Abstract.-}
The simplest patterns of qualitative changes on the {\it
configu\-ra\-tions of lines of principal curvature} around umbilic
points on surfaces whose immersions into $\mathbb R^3$ depend
smoothly on a real parameter (codimension one umbilic
bifurcations) are described in this paper. Global effects, due to
umbilic bifurcations, on these configurations such as the
appearance and annihilation of periodic principal lines, called
also principal cycles,  are also studied here.}

\vskip .4cm

\section{Introduction}\label{sec:intro}

 Let $\mathbb M^2$ be a connected, compact, oriented, two dimensional smooth
manifold  and let $\mathbb R^3$ be the {\it Euclidean} $\, 3-$ {\it space},  endowed
with a once for all fixed {\it orientation}, {\it Euclidean  inner product}: $<\, ,\, >$
 and {\it norm}: $|\,\;\,|  =  <\, ,\, >^{1/2}$.

An immersion  $\alpha$ of $\mathbb M^2$ into  $\mathbb R^3$ is  a map such that the
derivative  $D{\alpha}_p $ from  ${\mathbb T}{\mathbb M}^2_p$  to ${\mathbb R}^3$ is one
to one, for every
 $p \in {\mathbb M}^2.$
 Here ${\mathbb T}{\mathbb M}^2_p$ stands for the tangent space  of ${\mathbb M}^2$ at $p$.
  The tangent bundle of $\mathbb M^2$ will be denoted by ${\mathbb T}{\mathbb M}^2$. The
  tangent projective bundle
will be denoted by ${\mathbb P}{\mathbb M}^2$ and its projection onto $ {\mathbb M}^2$
will be designated by $\Pi$. The fiber ${\mathbb P}{\mathbb M}^2_p$ = ${\Pi^{-1}}(p)$
is the projective line over $p$  which is the space of tangent lines through $p$.

 Denote by  $\mathcal I^r$= ${\mathcal I}^r({\mathbb M}^2,{\mathbb R}^3)$ the set
  of $C^r$-immersions of ${\mathbb M}^2$ into ${\mathbb R}^3$.  When endowed with
   the $C^s,\, s\leq r,\,$ topology this space will be denoted
   by $\mathcal I^{r,s}$= ${\mathcal I}^{r,s}({\mathbb M}^2,{\mathbb R}^3).$

 To every $\alpha \in {\mathcal I}^r$ is associated its {\it Gaussian normal
  map} $N_\alpha : {\mathbb M}^2 \to {\mathbb S}^2$, with values in the
  unit 2-sphere ${\mathbb S}^2$,
defined  by
$$N_{\alpha}(p) = (|{\alpha}_u\wedge{\alpha}_v|^{-1}) \alpha_u \wedge \alpha_v , $$
\noindent  where  $(u,v): ({\mathbb M}^2, p) \to  ({\mathbb R}^2
,0)$ is a local positive chart on an open set of ${\mathbb M}^2$
containing  $p$; $ \wedge $  denotes the {\it vector}  or {\it
wedge product} in the oriented space ${\mathbb R}^3$; ${\alpha}_u
= \partial{\alpha}/\partial u$ and
${\alpha}_v=\partial{\alpha}/\partial v$. Clearly  $N_{\alpha}$ is
well defined and of class $C^{r-1}$  in  ${\mathbb M}^2$.

Since  $DN_\alpha (p)$  has its image contained in the image of
$D\alpha(p)$,  the {\it Weingarten endomorphism} $\omega_{\alpha}:
{\mathbb T}{\mathbb M}^2 \to  {\mathbb T}{\mathbb M}^2$
  is well defined by  $$(D\alpha)\,\omega_\alpha  =  DN_\alpha.$$

  It is
well known (\cite {Sp, St})  that  $\omega_\alpha$ is  self adjoint,
 when $ {\mathbb T}{\mathbb M}^2$  is endowed with the
 {\it First Fundamental Form}, given by

 $$I_\alpha (\, . , .\,) = <D\alpha (\, . \,), \,D\alpha (\, . \,)>.$$

Let ${\mathcal K}_\alpha$= $\det(-\omega_\alpha)$  and
${\mathcal H}_\alpha$= $\frac{1}{2}${\it trace}$(-\omega_\alpha)$  be respectively the
{\it Gaussian} and  {\it Mean Curvatures}  of the immersion  $\alpha$.

If $({\mathcal H}_\alpha^2-
{\mathcal K}_\alpha )(p)=0$, the point  $p \in {\mathbb M}^2$ is
called an {\it umbilic point} of $\alpha$.  The (closed) set of
 umbilic points of  $\alpha$ will be denoted by
${\mathcal U}_\alpha$.  The eigenvalues $k_{1,\alpha},
 \,k_{2,\alpha}$ of $-\omega_\alpha$ are always real and  given by

$$k_{1,\alpha} = {\mathcal H}_\alpha -({\mathcal H}_\alpha^2 -
{\mathcal K}_\alpha)^{1/2}, \; \;\;     k_{2,\alpha} =
{\mathcal H}_\alpha +({\mathcal H}_\alpha^2 - {\mathcal K}_\alpha)^{1/2}.$$

They are  called respectively the  {\it minimal} and
 {\it maximal principal curvatures} of  $\alpha$.
  It holds that $k_{1,\alpha} < k_{2,\alpha}$, except
    on ${\mathcal U}_\alpha$, where $k_{1,\alpha}\, =\,k_{2,\alpha}$.

The {\it Second Fundamental
Form} of $\alpha$ of is defined by

 $$ II_{\alpha}(\,. \, , \, .\, ) = - <\omega_{\alpha}(\, .\,),\, . \, >.  \eqno$$

The {\it Normal Curvature},  $k_{n,\alpha}(p,{\it l})$,  on a line $\it l$
through the point $p$, is defined by

 $$ k_{n,\alpha}(p,{\it{l}})= \frac{II_{\alpha}(\,. \, , \, .\, )}{I_\alpha (\, . , .\,) }, \eqno$$

\noindent evaluated at any  vector generating the line  ${\it {l}}$.

It is well known (\cite {Sp, St}) that  at each point  $p \in {\mathbb M}^2$, $k_{1,\alpha}(p)$ is the
minimum  and  $k_{2,\alpha}(p)$ is the  maximum of  $k_{n,\alpha}(p,{\it{l}})$, taken on all
 tangent lines $\it l$ through the point $p$.

The eigenspaces of  $-\omega_\alpha$ associated to the principal curvatures
define on
  ${\mathbb M}^2  \setminus {\mathcal U}_\alpha$ two  $C^{r-2}$  line fields  ${\mathcal L}_{1,
\alpha}$ and  ${\mathcal L}_{2,
\alpha}$   called
respectively the {\it  minimal} and {\it maximal principal line fields} of  $\alpha$, which  are
mutually orthogonal in  ${\mathbb T}{\mathbb M}^2$,   with the metric  $I_\alpha$.

The {\it principal line fields} are characterized by Rodrigues'
equations \cite {St, Sp}:

$${\mathcal L}_{i,
\alpha} = \{v\in {\mathbb T}{\mathbb M}^2;\; \omega_{\alpha}(v) + k_{i,\alpha}v = 0\}\;\; i=1,\,2.\eqno$$

Elimination of $k_{i,\alpha},\;\; i=1,\,2$ in these  equations leads to a single
quadratic differential equation $\tau_{g,\alpha}=0$ for the principal lines of $\alpha$,
in terms of  the {\it geodesic torsion}  of $\alpha$, defined  in the direction $(.)$ as

$$\tau_{g,\alpha}(.,.)=<DN_{\alpha}(.)\wedge D{\alpha}(.),N_{\alpha}>. \eqno$$

 Calculation (see \cite{Sp, St}) shows that in a  chart
$(u,v)$  in   ${\mathbb M}^2$,  $\tau_{g,\alpha}$
has the form

\begin{equation} \label{eq:1d}
\tau_{g,\alpha}([du:dv])= \frac{(Fg - G f)dv^2 +(Eg-Ge)dudv +
(Ef-Fe)du^2}{(E du^2+2Fdudv+Gdv^2)\,\sqrt{EG-F^2}}.
   \end {equation}

The coefficients are functions  of $(u,v)$, characterized in terms of the  {\it first}
and {\it second fundamental forms} of the immersion $\alpha$, written in the chart
$(u,v)$ as follows:
$$I_\alpha= <D\alpha,D\alpha>= E du^2+2Fdudv+Gdv^2,$$
\noindent  with $E=<\alpha_u,\alpha_u>$, $F=<\alpha_u,\alpha_v>$, $G=<\alpha_v,\alpha_v>$, and
$$II_{\alpha}= -<DN_{\alpha},\,D\alpha>=<N_{\alpha} , D^2\alpha> = e du^2 + 2f dudv + g dv^2,$$
\noindent  with $e=<N,\alpha_{uu}>=-<N_u,\alpha_u>$,
$f=<N,\alpha_{uv}>=-<N_u,\alpha_v>$, $g=<N,\alpha_{vv}>=-<N_v,\alpha_v>$, where $N=N_{\alpha}$.

\vskip 0.2cm

    The integral curves of   ${\mathcal L}_{i,
\alpha}\, , \; i=1,\, 2$  are called {\it lines of
 minimal}, $i=1$, and  {\it maximal}, $i=2$,  {\it principal curvature} of $\alpha$.  The family of such
curves i.e. the {\it integral foliation} of  ${\mathcal L}_{i,
\alpha}, \; i=1,\, 2,$  on   $\mathbb M^2\setminus{\mathcal U}_\alpha$ will be
denoted by  ${\mathcal F}_{i,\alpha}, \; i=1,\, 2,$  and will be called the {\it
 minimal}, $i=1$, and  {\it maximal}, $i=2$,   {\it principal foliations} of $\alpha$.

    By the {\it principal
configuration  of}  $\alpha$ will be meant  the following triple
$${\mathcal P}_\alpha= ({\mathcal U}_\alpha, {\mathcal F}_{1,
\alpha}, {\mathcal F}_{2,
\alpha}).$$

    The local study of principal configurations around an umbilic
point  received considerable attention in the classical works of  Monge \cite{Mg}, Cayley \cite{Ca},
Darboux \cite{Da} and Gullstrand \cite{Gu}, among others.

More recently Palmeira \cite {Pa}, Gui\~nez and Gutierrez
\cite{gugu}, \cite{GG},  Gui\~nez \cite{guin},  \cite {gui}, Bruce
 and Fidal \cite {B-F} and  Bruce and Tari \cite {B-T}, S\'anchez-Bringas and Galarza \cite{sb}, Mello \cite{me},
 to mention just a few,  have
studied solutions of more general binary differential equations
near singularities.

An initial step to study the generic properties of principal
configurations on algebraic surfaces has been given  in \cite
{Ga-S1}.

    The study of the global features of principal configurations  ${\mathcal P}_\alpha$
which remain topologically undisturbed under small perturbations of the
immersion  $\alpha$ --{\it principal structural stability}-- was initiated  by Gutierrez
 and Sotomayor in \cite{G-S1, G-S2, G-S3}. There were established sufficient
conditions for immersions  $\alpha$ of class  $C^r$  of a compact oriented
surface ${\mathbb M}^2$ into  $\mathbb R^3$ to have a $C^s$-{\it structurally  stable
 principal
configuration}, $ r > s \geq 3$.  This means that for any immersion    $\beta$
sufficiently $C^s$-close to  $\alpha$,  there must exist a homeomorphism  $h$  of
${\mathbb M}^2$  which maps  ${\mathcal U}_ \alpha$ onto  ${\mathcal U}_ \beta$   and
 maps lines of ${\mathcal F}_{i,\alpha}, \; i=1,\, 2,$ to  those of   ${\mathcal F}_{i,
\beta}, \; i=1,\, 2$.

    The conditions for principal structural stability established in \cite  {G-S1, G-S2, G-S3}
     are reviewed in  sections \ref{sec:3} and \ref{sec:4}. They
are imposed on each of the topological invariant entities of the principal
configuration.  Namely,

$a$) the umbilic points,

$b$) the periodic lines of curvature, also called principal cycles,

$c$) the absence  of umbilic separatrix connexions,  and

$d$) the absence of non-trivially  recurrent principal lines.

\vskip 0.2cm

    For principal structural stability, the umbilic points, which are
regarded as the singularities of the principal configuration, are
assumed to be {\it Darbouxian} \cite {Da, G-S1, G-S3}.

\begin{figure}[htbp]
 \begin{center}
 \hskip .5cm
  \includegraphics[angle=0, width=12cm]{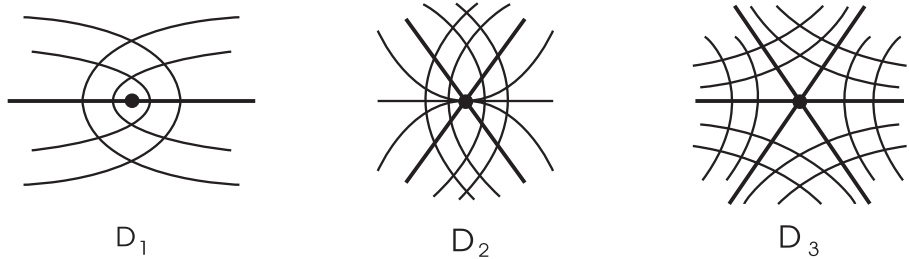}
  \caption{Principal curvature lines near the umbilic points $D_i$
 and their separatrices \label{fig:darboux}}
 \end{center}
  \end{figure}

In Section \ref{sec:2}, the meaning of this assumption is stated
intrinsically  and in terms of local coordinates involving the
coordinate expression for the third order jet of the immersion at
an umbilic point. See  Fig. \ref{fig:darboux} or an illustration
of
 the principal configuration around a Darbouxian umbilic. The subscript
  stands for the number of umbilic separatrices approaching the point. This
   number is the same for both the minimal and maximal principal curvature foliations.

\vspace{-.3cm}

    The purpose  of this paper is to study the simplest qualitative
changes --bifurcations-- exhibited by  the principal
configurations under small perturbations of an immersion which
violates in a minimal fashion the  Darbouxian structural stability
condition on umbilic points.

 This leads to two types of
umbilic patterns: $D^{1}_{2}$  and $D^{1}_{2,3}$, illustrated in
Fig. 2. The superscript stands for the codimension which is the
minimal  number of parameters on which depend the families of
immersions exhibiting
 persistently the pattern.    The subscripts stand for the number of
  separatrices approaching the umbilic. In the first case, this number
   is the same for both the minimal and maximal principal curvature
   foliations. In the second case, they are not equal and, in our notation,
   appear  separated by a comma.

 \begin{figure}[htbp]
 \begin{center}
 \hskip .5cm
  \includegraphics[angle=0, width=12cm]{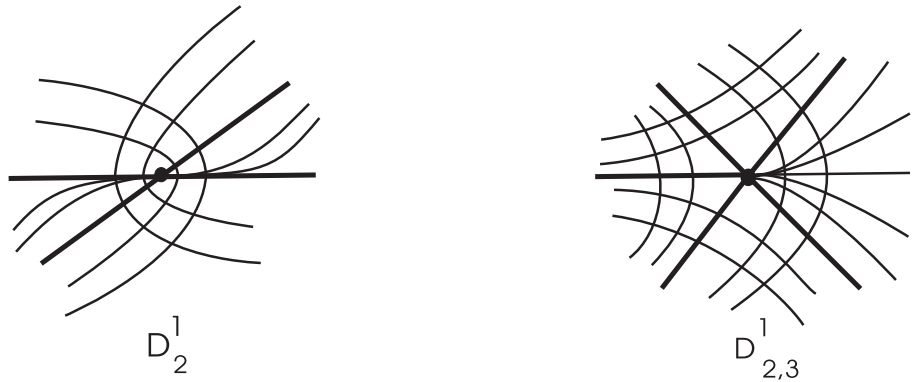}
  \caption{\label{fig:semidarboux} Principal curvature lines near the umbilic points  $D^{1}_{2}$, left, and
  $D^1_{2,3}$, right,
 and their separatrices. }
 \end{center}
  \end{figure}
\vspace{-.3cm}

The precise definitions and  bifurcation analysis of these points
is carried out in a local context in Section \ref{sec:2}.  The
global implications on the bifurcations  of principal cycles are
studied in sections 3 and 4. This is preceded by a short review of
the stability conditions, items   $a)$, $b)$,  $c)$ and $d)$,
since they are essential to state properly the above mentioned
global implications  on the appearance and annihilation of
principal cycles as well as for the statements concerning the {\it
relative principal stability } of the  immersions at umbilic
bifurcations, also included in this paper.  This is also done in
sections 3 and 4. Section \ref{cr} links the results of this paper with
other bifurcations of principal configurations discussed in \cite{G-S7}.

The  study presented here  is motivated by   the {\it Theory of
First Order Structural Stability} of vector fields on surfaces due
to Andronov-Leontovich \cite{A-L} and further developed by
Sotomayor \cite{S}. The conditions $a)$, $b)$,  $c)$ and $d)$,
above, are the counterpart, for  the  structural stability of
principal configurations,  of those  of Andronov-Pontrjagin \cite
{A-L} and Peixoto \cite {pe}, for the structural stability of
vector fields on surfaces. Particularly  meaningful here is the
openness and density theorem due to Peixoto \cite{pe}.

The analogy is pursued here to the case of {\it first order
structural stability}, following the steps of the work of
Sotomayor \cite{S}, where the transversality to key Banach
submanifolds are instrumental for the geometric explanation of the
codimension one generic bifurcation patterns. Elaboration of these
ideas  gives rise in section \ref{sec:4} to the definition of
certain Banach submanifolds consisting  of immersions where the
simplest umbilic bifurcations happen.

The Theory of Bifurcations of general vector fields has been
recently developed in several promising  directions. See the books
of   Chicone \cite{chi},  and Chow and Hale \cite{ch},
Guckenheimer and Holmes \cite{gho},  IIyashenko and Li W. Gu \cite
{I-L} and Roussarie \cite{rr}, to mention just a few.

    The bifurcations of Principal Configurations due to the minimal
 violation
of the structural stability conditions  imposed on principal
cycles and umbilic separatrix connections have been studied in
\cite{G-S4, G-S5, Ga-S2, G-S7}. To the knowledge of the authors, a
systematic  analysis of the violation of condition  $d$ --leading
to non-trivial  recurrences of principal lines-- has not yet been
carried out. Explicit examples of such  recurrences -- missing in
the classical geometry  literature -- have been provided in \cite{G-S2,
G-S3}. Their classification, however, seems quite distant.

\section { Umbilic Points and their Generic
Bifurcations}\label{sec:2}

\subsection{ Preliminaries concerning umbilic points} \label{ss:2.1}

Denote by ${\mathbb P}{\mathbb M}^2$ the projective tangent bundle over ${\mathbb M}^2$,
with projection ${\Pi}$. For any chart $(u,v)$ on an
open set U of ${\mathbb M}^2$ there are defined two charts $(u,v;p=dv/du)$ and
$(u,v;q=du/dv)$ which cover ${\Pi}^{-1}(U)$.

The equation \ref {eq:1d} of principal lines, being quadratic, is
well defined in the projective bundle. Thus,
 for every $\alpha$ in ${\mathcal I}^r$,

 $${\mathbb L}_\alpha=\{ \tau_{g,\alpha}=0, \}$$

\noindent defines a variety
 on ${\mathbb P}{\mathbb M}^2$, which is regular and of
class $C^{r-2}$  over ${\mathbb M}^2 \setminus {\mathcal U}_\alpha$. It doubly
covers ${\mathbb M}^2 \setminus {\mathcal U}_\alpha$ and contains a projective
line ${\Pi}^{-1}(p)$ over each point $p \in {\mathcal U}_\alpha$.

\begin{definition}\label{df:da}
A point $p \in {\mathcal U}_\alpha$ is {\it Darbouxian} if  the following two conditions hold:
\begin{itemize}
\item [$T:$] The variety  ${\mathbb L}_\alpha$ is regular also
over ${\Pi}^{-1}(p)$. In other words,  the derivative of
$\tau_{g,\alpha}$ does not vanish on the points of projective line
${\Pi}^{-1}(p)$. This means that the derivative in directions
transversal to ${\Pi}^{-1}(p)$ must not vanish. \item [$D:$] The
principal line fields ${\mathcal L}_{i,\alpha}, \,i=1,2$ lift to
a single line field ${\mathcal L}_{\alpha}$ of class $C^{r-3}$,
tangent to ${\mathbb L}_\alpha$, which extends to a unique one
along  ${\Pi}^{-1}(p)$, and there it has only hyperbolic
singularities, which must be either
 \begin{itemize}
\item [$D_1:$] a unique saddle \item [$D_2:$] a unique node
between two saddles, or
 \item [$D_3:$] three  saddles.
\end{itemize}
\end{itemize}
\end{definition}

    For calculations it will be essential  to  express the Darbouxian conditions in  a
     Monge    local chart $(u,v)$: $({\mathbb M}^2, p) \to  ({\mathbb R}^2,0)$
     on $\mathbb M^2$ , $p\in {\mathcal U}_\alpha$,
as  follows.

Take an  isometry  $\Gamma$
of  $\mathbb R ^3$ with $\Gamma(\alpha (p)) = 0$  such that $\Gamma(\alpha (u,v))=(u,v,h(u,v))$, with

\begin{equation}\label{eq:2.1}
 \aligned
 h(u,v) =& \frac k2 (u^2+v^2 ) +  (a/6)u^3 + (b/2)uv^2 + (b'/2) u^2v\\
 +&(c/6)v^3 +(A/24)u^4 + (B/6)u^3v  \\
+&(C/4)u^2 v^2 + (D/6)uv^3 + (E/24)v^4  + O((u^2 +v^2)^{5/2}).
\endaligned
\end{equation}

    To obtain simpler  expressions assume  that the coefficient $ b'$ vanishes.

 This
is achieved by means of a suitable  rotation in the $(u,v)$-plane.

    In the affine chart $(u,v;\, p=dv/du)$ on ${\mathbb P}({\mathbb M}^2)$ around ${\Pi}^{-1} (p)$, the
variety  $ {\mathbb L}_\alpha $ is given by the following equation.

\begin{equation}\label{eq:2.2}
{\mathcal T}(u,v,p)=  L(u,v)p^2   + M(u,v)p + N(u,v) = 0, \,p=dv/du.
\end{equation}

 According to \cite{Sp, St},  the functions $L$, $M$ and $N$ are obtained   from
  equations \ref{eq:1d} and \ref{eq:2.1} as follows:
$$\aligned L &= h_uh_v h_{vv}- (1+h_v^2)h_{uv}\\
M &= (1+h_u^2)h_{vv}- (1+h_v^2)h_{uu}\\
N &= (1+h_u^2)h_{uv} - h_u h_v h_{uu}.\endaligned $$

Calculation taking into account the coefficients in equation \ref {eq:2.1}, with $b' =0$, gives:

\begin{equation}\label{eq:2.2.1}
\aligned
L(u,v) =&-bv - B/2(u^2) -  (C-k^3)uv - D/2(v^2) + M_1^3(u,v) \\
M(u,v)=&(b-a)u+ cv+ [(C-A)/2+k^3]u^2+ (D-B)uv  \\
+& [(E-C)/2-k^3]v^2 +M_2^3(u,v)\\
N(u,v) =& bv + B/2u^2 + (C-k^3)uv +  D/2v^2 + M_3^3(u,v),
\endaligned
\end{equation}

\noindent with $M_i^3(u,v) = O((u^2+v^2)^{3/2})$,   i = 1, 2, 3.

    These expressions are obtained
from the calculation of the coefficients of the first and second
fundamental forms in the chart $(u,v)$ and substitution into
\ref{eq:1d}. See also \cite {Da,  G-S1, G-S3}.  With longer
calculations,  Darboux \cite {Da} gives the full expressions for
any value of  $b'$.

\begin{remark} \label{rm:dar}
     The regularity condition $T$ in definition \ref{df:da} is equivalent to impose  that
        $b(b-a)\neq 0$.  In fact, this inequality also implies regularity at $p=\infty$. This
can be seen in the chart $(u,v;\,q=du/dv)$, at $q=0$.

Also this condition is equivalent to the transversality of the
curves $M=0, \, N=0$
 \end{remark}

    The line field ${\mathcal L}_\alpha$ is expressed in the chart $(u,v;\,p)$  as being generated by the vector
field $X = X_\alpha$, called the {\it Lie-Cartan}  vector field of
equation \ref{eq:1d}, which is tangent to ${\mathbb L}_\alpha$ and
is given by:

\begin{equation} \label {eq:lc}
\aligned
             \dot u=& \partial {\mathcal T}/\partial p\\
   \dot v= & p \partial {\mathcal T}/\partial p\\
           \dot  p =&-({\partial {\mathcal T}}/{\partial u} + p{\partial {\mathcal T}}/{\partial v})
\endaligned
\end{equation}

Similar expressions hold for the chart $(u,v;q=du/dv)$ and the pertinent vector field  $Y = Y_\alpha$.

The function ${\mathcal T}$ is a first integral of $X = X_\alpha$.
The projections of the integral curves of $X_{\alpha} $ by
$\Pi(u,v,p)=(u,v)$ are the lines of curvature.
 The singularities of $X_{\alpha} $ are given by $(0,0,p_i)$ where $p_i$ is a root of the
  equation $p(bp^2-cp+a-2b)=0$.

    Assume that $b\neq 0$, which occurs under the regularity condition $T$, then the
    singularities of $X_\alpha$ on the surface ${\mathbb L}_\alpha$ are located on
the $p$-axis at the points with coordinates $p_0,\,p_1,\, p_2$

\begin{equation} \label{eq:lcs}
\aligned
    p_0=&0,\\
    p_1=&c/2b - \sqrt{ (c/2b)^2 - (a/b) + 2},\\
    p_2=&c/2b + \sqrt{ (c/2b)^2 - (a/b) + 2} \endaligned
\end{equation}

\begin{remark} \label{rm:dard} \cite{G-S1}
  Assume  the notation  established in equation \ref{eq:2.1}.
Suppose that the transversality condition $T:  b(b-a)\ne 0$ of definition \ref{df:da} and
remark \ref{rm:dar} holds. Let
$ \Delta =-[ (c/2b)^2 - (a/b) + 2].  $
Calculation of the hyperbolicity conditions for singularities \ref{eq:lcs} of the vector
 field \ref{eq:lc} --see \cite{G-S1}-- have led  to establish the   following equivalences:
\begin{itemize}
\item[$D_1$)]   $\equiv \;\;\Delta>0$

\item[$D_2$)]$\equiv \;\;\Delta<0 $ and $ \;\; 1< \dfrac ab  \ne
2$

\item[$D_3$)]    $\equiv \;\; \dfrac ab < 1. $
\end{itemize}
\end{remark}

\begin{figure}[htbp]
 \begin{center}
 \includegraphics[angle=0, width=11cm]{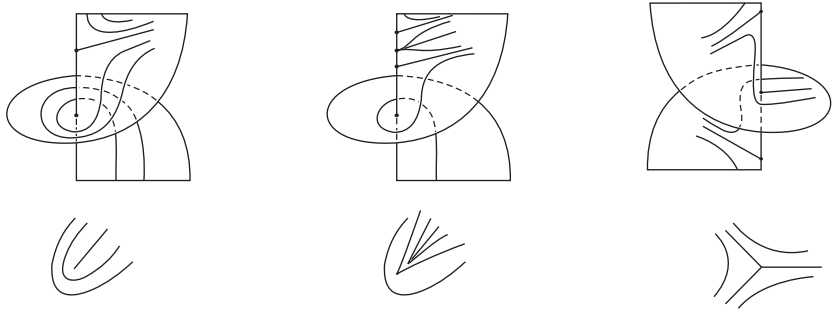}
 \caption{Darbouxian Umbilic Points, corresponding ${\mathbb L}_\alpha$ surface and lifted
  line fields.\label{fig:darbouxp}}
   \end{center}
 \end{figure}

    See Figs. \ref{fig:darboux} and \ref{fig:darbouxp} for an illustration of the three possible types of
Darbouxian umbilics. The distinction between them is expressed in terms of the coefficients
 of the  $3$-jet of equation \ref{eq:2.1}, as well as in the lifting
of singularities to the surface ${\mathbb L}_\alpha$. See remarks
\ref{rm:dar} and \ref{rm:dard}.

    The subscript  $i = 1,2,3$  of  $D_i$  denotes the number of {\it umbilic separatrices}
    of $ p$.
 These are principal lines which tend to the umbilic
point  $p$  and  separate regions of different patterns of approach to it.
For Darbouxian points, the umbilic separatrices are the projection
into $\mathbb M^2$ of the saddle separatrices transversal to the projective line over
the umbilic point.

 It can be proved that the only umbilic points for which  $\alpha \in  {\mathcal I}^r$
is locally $C^s$-structurally stable, $ r > s \geq 3,$ are the Darbouxian ones. See  \cite{B-F, G-S3}.

\vskip .4cm

 The implicit surface ${\mathcal T} (u,v,p)=0$ is regular in a neighborhood of the
 projective line  if and only if $b(b-a)\ne 0$. Near the singular point $p_0=(0,0,0)$
  of $X_{\alpha}$ it follows that ${\mathcal T}_v(p_0)=b\ne 0$ and therefore, by the
  Implicit Function Theorem, there exists a function $v$ such that   ${\mathcal T}(u,v(u,p),p)=0$.
   The function $v=v(u,p)$   has the following Taylor expansion
$$v(u,p)=-\frac {B}{2b} u^2 +\frac{a-b}b up + O(3).$$

For future reference we record the expression the vector field $X_{\alpha}$
 in the chart $(u,p)$.

 \begin{equation} \label{eq:lci}
\aligned \dot u =& {\mathcal T}_p(u,v(u,p),p)\\
=&(b-a)u+\frac 12\frac{[b(C-A+2k^3)-cB]}{b} u^2 +\frac{ c (a-b)}{b} u p +O(3) \\
\dot p =&-({\mathcal T}_u+p{\mathcal T}_v)(u,v(u,p),p) = -Bu+(a-2b)p   -cp^2  \\
+&\frac 12 \frac{  B (C-k^3)-a_{41}b}{b} u^2 +\frac{ b(A-C-2k^3)+a(k^3-C)}{b} up+ O(3). \endaligned
\end{equation}
\noindent where $a_{41}$ is term of order five in $h$ which, however, will have no influence in what follows.
\vskip .3cm

    Two generic patterns of bifurcations of umbilic points appear in this work. The first
     one  occurs due  to the violation of the Darbouxian condition
$D$, while $T$ is preserved, leading to the pattern called $D^{1}_{2}$, studied in subsection \ref{ss:2.2}. The
 second one  happens due to the violation  of condition $T$, leading to
the  pattern  denominated $D^{1}_{2,3}$, studied in subsection \ref{ss:2.3}.

 The meaning of the above
mentioned  genericity assertion for these bifurcations  is made precise  in subsection  \ref{ss:2.4},
 Theorem \ref{th:1}.


\subsection{The $D_{2}^1$ Umbilic Bifurcation Pattern}\label{ss:2.2}

Here will be studied the qualitative changes - bifurcations - of
the principal configurations around non Darbouxian umbilic points
such that the regularity (or transversality) condition  $T: b(a-b)
\neq  0$,  which implies their isolatedness, is preserved  and
only the condition $D$  is violated in the mildest possible way.

\begin{definition}\label{df:da12}
A point $p \in {\mathcal U}_\alpha$ is said to be of type $D^{1}_{2}$ if  the following holds:
\begin{itemize}
\item [$T:$] The variety  ${\mathbb L}_\alpha$ is regular also
over ${\Pi}^{-1}(p)$. In other words,  the derivative of
$\tau_{g,\alpha}$ does not vanish on the points of projective line
${\Pi}^{-1}(p)$. This means that the derivative in directions
transversal to ${\Pi}^{-1}(p)$ must not vanish. \item [$D_{2}^1
:$] The principal line fields ${\mathcal L}_{i,\alpha}, \,i=1,2$
lift to  a single line field ${\mathcal L}_{\alpha}$ of class
$C^{r-3}$, tangent to ${\mathbb L}_\alpha$, which extends to a
unique one along  ${\Pi}^{-1}(p)$, and there it has a hyperbolic
saddle  singularity and a saddle-node whose central manifold is
located along the projective line over $p$.
\end{itemize}
\end{definition}

In coordinates $(u,v)$, as in the notation above,
this means that

$T\,:$ $b(a-b) > 0$ and  either

 $1)$ $a/b = (c/2b)^2+2$,  or

$2)\,$ $a/b = 2$.

We point out that due to the particular representation of the {\it 3-jets} taken here,
 with $b'=0$, the space $a,b,c$ in the case $2)$ is not transversal, but tangent,
  to the  manifold of jets with $D^{1}_{2}$ umbilics.

\begin{remark} \label {rm:sep}
The $D^{1}_{2}$  umbilic point has two separatrices.

The  {\it isolated} one is characterized by the fact that no other principal
line which approaches the umbilic point is tangent to it. 

The other separatrix, called {\it non-isolated}, has the property
that every principal line distinct from the isolated one, that
approaches the point does so tangent to it.

\end{remark}

These separatrices bound the parabolic sector of lines of curvature
approaching the point; they also constitute the boundary of the
hyperbolic sector of the umbilic point.

\begin{figure}[htbp]
 \begin{center}
 \includegraphics[angle=0, width=12cm]{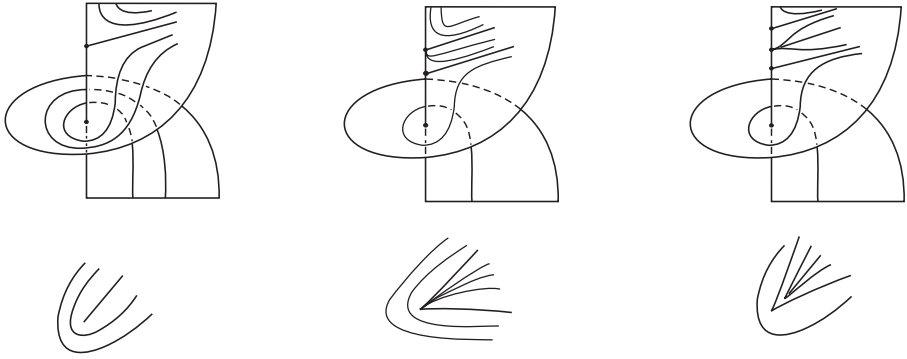}
 \caption{Umbilic Point $D^{1}_{2}$ , corresponding  ${\mathbb L}_\alpha$ surface
 and lifted line fields.\label{fig:d12bifu}}
  \end{center}
 \end{figure}

    The bifurcation illustrated in Fig. \ref{fig:d12bifu} shows that the non-isolated
separatrix disappears when the point  $D^{1}_{2}$  changes to  $D_1$ and that it
turns into an isolated  $D_2$  separatrix when it changes into  $D_2$ .
 It can be said that  $D^{1}_{2}$  represents the simplest transition
between $D_1$  and  $D_2$   Darbouxian umbilic points, which
occurs through the annihilation of an umbilic separatrix -- the
non-isolated one --.

With the notation in equations \ref{eq:2.2} and \ref{eq:2.2.1}, write

\begin{equation}
\aligned
     L=& -bv+M^2_1(u,v), \\
     M=&(b-a)u + cv +M^2_2(u,v),\\
     N =&  bv+M^2_3(u,v),
\endaligned
\end{equation}

\noindent with $M^2_i(u,v)= O(u^2+v^2)$.

 Condition $D^{1}_{2}$ is equivalent to the existence of a non zero double root for
$bp^2-cp+a-2b=0$, which amounts to $b\neq 0$ and $p_1 = p_2\neq p_0$.

    Assuming $b(b-a)\neq 0$, the curves $L=0$ and $M =0$ meet
transversally at $(0,0)$ if and only if $b\neq a$. It was shown in \cite{G-S1} that $D_1$
is satisfied if and only if the roots of $ bp^2-cp+a-2b = 0$ are non vanishing
and purely imaginary.

Also, $D_2$ is satisfied if and only if $bt^2-ct+a-2b=0$
has two distinct non zero real roots, $p_1,\, p_2$  which verify $p_1 p_2  > -1$.

This means that the rays  tangent to the separatrices
are pairwise distinct and contained
in an open right angular sector.

The local configuration of $D^{1}_{2}$
is established now.

     \begin {proposition} \label {pr:2.2.1)}
 Suppose that $\alpha \in {\mathcal I}^r, r\geq 5$, satisfies condition $D^{1}_{2}$
 at an  umbilic point $ p$.
Then the local principal configuration of $\alpha$ around $ p$ is
that of Fig. \ref{fig:semidarboux}, right and Fig.
\ref{fig:d12bifu}, center.
\end{proposition}

\begin{proof}
Consider the Lie-Cartan lifting $X_{\alpha}$  as in equation
\ref{eq:lc},  which is of class $C^{r-3}$.  If $a=2b\ne 0$ and
$c\ne 0$, it follows that $p_0 =(0,0,0)$ is an isolated singular
 point of quadratic saddle node type with a center manifold contained in
 the projective line --the $p$ axis--. In fact, the eigenvalues of $DX_{\alpha}(0)$ are
 $\lambda_1=-b\ne 0$ and $\lambda_2=0$ and the  $p$ axis  is
 invariant; there
   $X_{\alpha}$, according to equation \ref{eq:lci}  is given by $\dot p=-cp^2+ o(2)$.

The other singular point of $X_{\alpha} $ is given by $p_1=(0,0,\frac cb)$.
It follows that
$$\aligned DX_\alpha(0,0,p_1)=\left[\begin{matrix} -b & -c & 0 \\ -c & -\frac{c^2}b &0 \\ A_1 & A_2 &\frac{c^2}b
\end{matrix}\right] \endaligned$$
where,
$$\aligned A_1=&\frac{ b^2c(A-k^3-2C)+bc^2(2B-D)+c^3(C-k^3)-b^3D }{b^3}\\
A_2=&\frac{  b^2c(B-2D)+bc^2(2C+k^3-E)+b^3(k^3-C)+Dc^3  }{b^3}\endaligned$$

The non zero eigenvalues of $DX(0,0,p_1)$ are
$\lambda_1=\frac{c^2}b, \;\;\; \lambda_2=-\frac{c^2+b^2}b$. In
fact,  $p_1$ is a hyperbolic saddle point of $X_{\alpha}$
having eigenvalues given by $\lambda_1$ and $\lambda_2$.

Similar analysis can be done when $(\frac{c}{2b})^2 - \frac ab +1 = 0$. In this case  $X_{\alpha} $
 and $p_1=(0,0,\frac cb)$ is a quadratic saddle node, with a local
 center manifold contained in the projective line. The point
$p_0=(0,0,0)$ is a hyperbolic saddle of $X_\alpha$. This case is
equivalent to the previous one, after a rotation in $(u,v)$ that sends de saddle-node to $p=0$.

\end{proof}

The principal configuration of the  $D^{1}_{2}$  umbilic point is
illustrated in Fig. 2, left,  and its bifurcation in Fig.
\ref{fig:d12bifu}.

\begin {proposition} \label {pr:2.2.2)}
 Suppose that $\alpha \in {\mathcal I}^r, r\geq 5$, satisfies condition $D^{1}_{2}$
 at an  umbilic point $ p$.
Then there is a function $\mathcal B$ of class $C^{r-3}$ on a
neighborhood $\mathcal V$ of $\alpha$ and a neighborhood $V$ of
$p$ such that  every $\beta \in \mathcal V $ has a unique umbilic
point $p_\beta$ in $V$.
\begin{enumerate}
\item[i)] $d{\mathcal B}(\alpha) \ne 0$

\item[ii)] ${\mathcal B}(\beta) >  0$ if and only if $p_\beta$ is
Darbouxian of type $D_1$

\item[iii)] ${\mathcal B}(\beta) <  0$ if and only if $p_\beta$
is Darbouxian of type $D_2$

\item[iv)] ${\mathcal B}(\beta) =  0$ if and only if $p_\beta$ is
of type $D^1_2$
\end{enumerate}

The principal configurations of $\beta$ around $ p$ is that of
Fig. \ref{fig:d12bifu}, left, right and center, respectively.
\end{proposition}

\begin{proof} Since $p$ is a transversal  umbilic point of $ {\alpha}$, the existence of the neighborhoods
 $\mathcal V$ and $V$
of $p_\beta$ follow from the Implicit Function Theorem. So we
assume that after   an isometry $\Gamma _\beta$ of $\mathbb R^3$,
with ${\Gamma_\beta} {\beta} (0)=0$,  in the  neighborhood $V$ are
defined coordinates $(u,v)$, also depending on $\beta$, on which
it is represented as:

$$h_{\beta}(u,v) = k_{\beta}/2(u^2+v^2) + a_{\beta}/6 (u^3) + b_{\beta}/2 (uv^2)
 + c_{\beta}/6(v^3) +   O({\beta}; (u^2+v^2)^4).$$

Define the function
$${\mathcal B}({\beta}) = [c_{\beta}/2 b_{\beta}]^2 - a_{\beta}/b_{\beta} +2,$$
\noindent whose zeros define locally the manifold of immersions with a $D^{1}_{2}$ point.

The derivative of this function in the direction of the coordinate $a$ is clearly non-zero.
\end{proof}

\subsection{The $D^{1}_{2,3}$ Umbilic Bifurcation Pattern} \label{ss:2.3}

    The second case of non-Darbouxian umbilic point  studied here, called $D^{1}_{2,3}$, happens when the regularity
condition
$T$ is violated.

\begin{definition}\label{df:d23}
An umbilic point  is said of type  $D^{1}_{2,3}$ if the transversality condition $T$
 fails at two points over  the umbilic point, at which
  ${\mathbb L}_\alpha$ is non-degenerate  of  {\it Morse } type.
\end {definition}

  \begin {proposition} \label {pr:2.2.1sn)}
 Suppose that $\alpha \in {\mathcal I}^r, \; r\geq 5$, and $p$ be  an umbilic
 point. Assume the notation in
\ref{eq:2.1} with $b=a\ne 0$ and $ b(C-A+2k^3)-cB\ne 0$.

Then  $ p$ is   of type  $D^{1}_{2,3}$ and the local principal
configuration of $\alpha$ around $ p$ is that of Fig. 2, right.

\end{proposition}

\begin{proof}  Consider the Lie-Cartan lifting $X_{\alpha}$  as in equation
\ref{eq:lc},  which is of class $C^{r-3}$.
 Imposing $a=b\neq 0$, by equations \ref{eq:lcs} and \ref{eq:lci},  the singular
 points of $X_{\alpha}$ are   $p_0$, $p_1$ and $p_2$, roots of the equation $p(bp^2-cp-b)=0.$

In fact, if $a=b\ne 0$,    it follows that $p_0$ is a quadratic saddle node  with  center
 manifold transversal to the projective line.

 From equation \ref{eq:lci} the eigenvalues are $\lambda_1=  0$ and $\lambda_2=-b$ and
the all center manifolds $W^c$ are tangent to the line $p=-\frac Bb u$. By invariant
 manifold theory it follows  that $X|W^c$ is local  topologically equivalent to

\begin{equation} \label{eq:xi}
\dot u=\frac 12\frac{[b(C-A+2k^3)-cB]}{b}
u^2+o(2):=-\frac{\chi}{2b}u^2+o(2).
\end{equation}

It follows that
$$\aligned DX_\alpha(0,0,p_i)=\left[\begin{matrix} 0 & -2bp_i+c & 0 \\ 0 & -p_i(2bp_i-c)
&0 \\ B_1 & B_2 &3bp_i^2-2cp_i-b \end{matrix}\right] \endaligned$$
\noindent where,
$$\aligned B_1=&(C-k^3)p_i^3 +(2B-D)p_i^2+(A-2C-k^3)p_i -B\\
B_2=&Dp_i^3 +(2C+k^3-E)p_i^2+(B-2D)p_i +k^3-C. \endaligned$$

The nonzero eigenvalues of $DX_\alpha(0,0,p_i)$ are
$\lambda_1=-2bp_i^2+cp_i=-b(p_i^2+1)$ and
$\lambda_2=3bp_i^2-2cp_i-b=b(p_i^2+1)$.

By invariant manifold theory, $(0,0,p_i)$ are saddles of
$X_{\alpha}$. The phase
 portrait of $X_{\alpha}$ near  these singularities are as shown in  Fig. \ref{fig:cone}.

\begin{figure}[htbp]
  \begin{center}
  \hskip 1cm
  \includegraphics[angle=0, width=8cm]{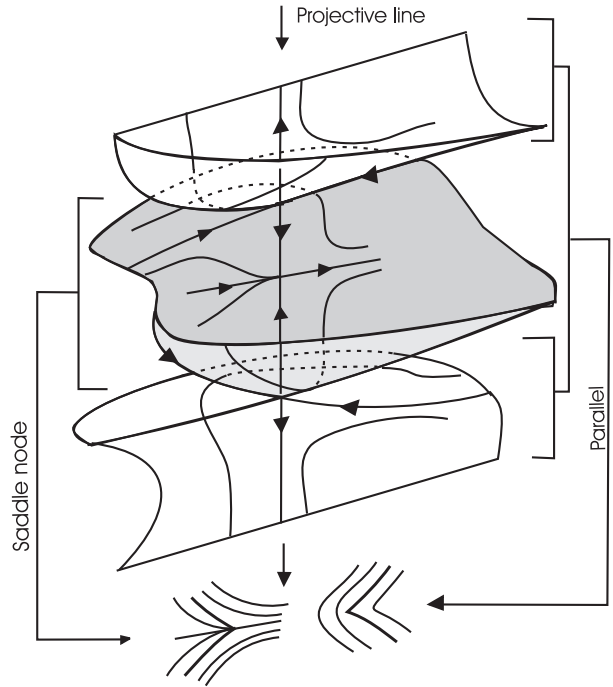}
  \caption{Lie-Cartan suspension $D^{1}_{2,3}$
  \label{fig:cone}}
    \end{center}
  \end{figure}

The two {\it critical} points $p_1$ and $p_2$ are  of conic type
on the variety ${\mathbb L}_\alpha$
 over  the umbilic point.

These points are  non-degenerate or of  {\it Morse } type,
according to the analysis below. At the points $(0,0,p_i)$ the
variety ${\mathcal T}(u,v,p)=0$ is not regular. In fact:

$\nabla{\mathcal T}(0,0,p)=[(b-a)p,-bp^2+cp+b,0].$

Therefore, for $a=b\ne 0$, at the two roots of the equation
$-bp^2+cp+b=0$ it follows
 that $\nabla{\mathcal T}(0,0,p_i)=(0,0,0), i = 1, 2.$

The Hessian of ${\mathcal T}$ at $(0,0,p_i)$ is

$$\aligned Hess({\mathcal T})(0,0,p_i) =\left[\begin{matrix}   \frac{p_i(-cB+b(C+2k^3-A))}{b}
& \frac{p_i( c(k^3-C)+b(D-B))}{b} & 0 \\
\frac{p_i( c(k^3-C)+b(D-B)}{b} & -\frac{ p_i(cD+b(C-E+2k^3))}{b}&-2bp_i+c\\
0 & -2bp_i+c & 0\end{matrix}\right]\endaligned$$

Direct calculation with the notation in equation \ref{eq:xi} gives,
 $$det (Hess({\mathcal T})(0,0,p_i))=
 \frac{p_i(-2bp_i+c)^2\chi}{b}=\frac{b}{p_i}(p_i^2+1)^2\chi\ne 0.$$

Therefore, $(0,0,p_i)$ is a non degenerate critical point of $\mathcal T$ of Morse
type and index $1$ or $2$ --a {\it cone}--, since ${\mathcal T}^{-1}(0)$ contains
 the projective line.
 \end{proof}

\begin{remark}\label{rm:d23eq}
Our analysis has shown  the   equivalence between the conditions
$a)$ and $b)$ that follow:

a) The non-vanishing on the
Hessian of ${\mathcal T}$ on  the critical points $p_1$ and $p_2$ over the umbilic.

b) The presence
of a {\it saddle-node}  at $p_0$ on the regular portion of the surface ${\mathbb L}_\alpha$, with
 central  separatrix  transversal to the projective
line over the umbilic.

Further direct calculation with equation \ref {eq:2.2.1} gives that  these two
conditions are equivalent to

c) The quadratic contact at the umbilic  between the curves $M=0$ and $N=0$.

\end{remark}
 In fact, from equation \ref{eq:2.2.1}  it follows that
 $M(u,v(u))=0$ for $v=-(B/2b) u^2+o(2)$ of class $C^{r-2}$.   Therefore $n(u)=N(u,v(u))$
 is of class $C^{r-2}$ and $ n(u)=  -(\chi/2b) u^2+o(2)$.

Notice also that, unlike the other umbilic points discussed here,
the two principal foliations around  $D^{1}_{2,3}$ are topologically distinct.

One
of them, located on the {\it parallel sheet},  has two umbilic separatrices and
 two hyperbolic sectors

The other, located on the {\it saddle-node  sheet},  has three umbilic
separatrices, one parabolic and two
hyperbolic  sectors.

The separatrix which is the common boundary of the hyperbolic
sectors will be called {\it hyperbolic separatrix}. See Figs.
\ref{fig:darbouxp}, \ref{fig:d12bifu}   and Fig.  \ref{fig:cone}
for illustrations.

    The bifurcation analysis describes the elimination
of two umbilic points  $D_2$ and $D_3$  which, under a deformation of the
immersion, collapse into a single umbilic point $D^{1}_{2,3}$ ,  and then, after a
further suitable arbitrarily small perturbation, the umbilic point is
annihilated.

\begin{figure}[htbp]
\begin{center}
 \includegraphics[angle=0, width=10cm]{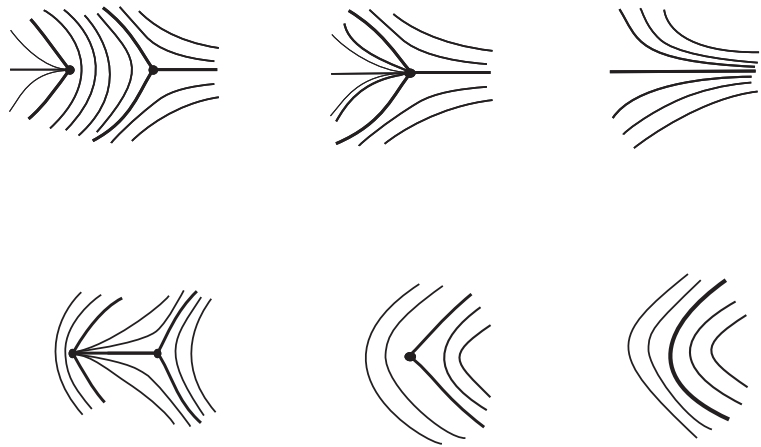}
 \caption{ Umbilic Point $D^{1}_{2,3}$  and bifurcation.\label{fig:bif23}}
  \end{center}
 \end{figure}

\begin {proposition}\label{pr:d23st}
 Suppose that $\alpha \in {\mathcal I}^r, r\geq 5$, satisfies condition $D^{1}_{2,3}$
 at an  umbilic point $ p$.
Then there is a function $\mathcal B$ of class $C^{r-3}$ on a
neighborhood $\mathcal V$ of $\alpha$ and a neighborhood $V$ of
$p$ such that

\begin{enumerate}
\item[i)]\;\; $d{\mathcal B}(\alpha) \ne 0$

\item[ii)]\;\; ${\mathcal B}(\beta) >  0$ if and only if $\beta$
has no umbilic points in $V$,

\item[iii)]\;\; ${\mathcal B}(\beta) <  0$ if and only if $\beta$
has two Darbouxian umbilic points of types  $D_2$ and $D_3$,

\item[iv)]\;\; ${\mathcal B}(\beta) =  0$ if and only if $\beta$
has only one   umbilic point in $V$, which is  of type
$D^1_{2,3}$.

\end{enumerate}

The principal configurations of $\beta$ around $ p$ are illustrated
in Fig. \ref{fig:bif23}, right, left and center, respectively.

\end{proposition}

\begin{proof}
Similar to that given in \cite{S}, page 15, for the saddle-node of vector fields,
 using the equivalence $c)$ of remark \ref{rm:d23eq}.  We define
 $\mathcal B$ as follows. An immersion $\beta $ in a neighborhood $\mathcal V$ of $\alpha$
 and a neighborhood $V$ of
  $p$ can be written in a Monge chart as a graph of a function
  $h_\beta(u,v)$. The umbilic points of $\beta$ are defined by the
  equation

 \begin{equation}\label{eq:ubeta} \aligned M_\beta &= (1+((h_\beta)_u)^2)(h_\beta)_{vv}-
 (1+((h_\beta )_v)^2)(h_\beta )_{uu}=0\\
                    N_\beta  &= (1+((h_\beta)_u)^2)(h_\beta )_{uv} -
                    (h_\beta)_u (h_\beta)_v (h_\beta)_{uu}=0.\endaligned
                    \end{equation}

For $\beta$ in a neighborhood of $\alpha$ it follows that
$M_\beta(u,v_\beta(u))=0$.

Define  ${\mathcal B}(\beta) = n_\beta(u_\beta)$,  where $u_\beta$ is the only
critical point of $n_\beta(u)  =N_\beta(u,  v_\beta (u))$.

Taking $h_\beta(u,v)=h(u,v)+\lambda uv$, where $h$ is as in
equation \ref{eq:2.1}   it follows
 by direct calculation that     $\frac{d{\mathcal
 B}(\beta)}{d\lambda}| _{\lambda=0}\ne 0$.
\end{proof}

The bifurcation of the point  $D^{1}_{2,3}$ can be regarded as the simplest
transition between umbilics $D_2$ and $D_3$ and non umbilic points.
See the illustration in Fig. \ref{fig:bif23}.

\subsection{A Transversality Theorem} \label{ss:2.4}

\begin {theorem} \label{th:1}
 In the space of smooth  mappings of ${\mathbb M}^2 \times {\mathbb R}
\to \mathbb R^3$ which are
immersions relative to the first variable,
those which have their umbilic
points Darbouxian, $D^1_2$ and $D^1_{2,3}$, forming a curve in
${\mathbb M}^2 \times {\mathbb R}$ whose
projection into $\mathbb R$  has only non-degenerate critical
 points at $D^{1}_{2,3}$, is open
and dense.
\end{theorem}

\begin{figure}[htbp]
  \begin{center}
  \hskip 1cm
  \includegraphics[angle=0, width=7cm]{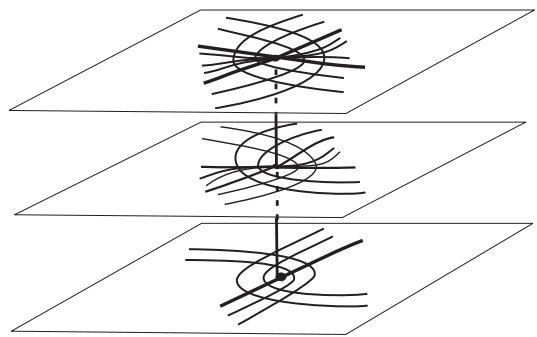}
  \caption{Bifurcation of   $D^{1}_{2}$
  \label{fig:d12sur}}
    \end{center}
  \end{figure}

\begin{figure}[htbp]
  \begin{center}
  \hskip 1cm
  \includegraphics[angle=0, width=13cm]{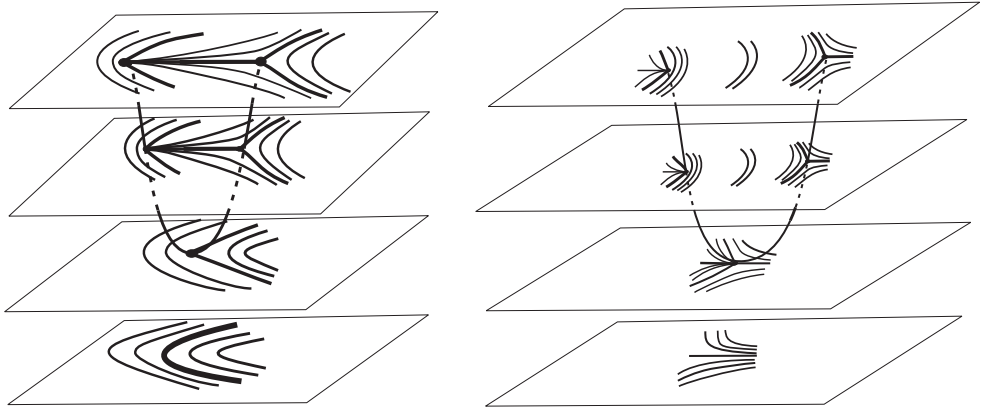}
  \caption{Bifurcation of   $D^{1}_{2,3}$
  \label{fig:d23sur}}
    \end{center}
  \end{figure}

\begin{proof}
Follows from Thom Transversality Theorem applied to the submanifold of three
jets of immersions at umbilic points, stratified by the Darbouxian,
having codimension 2, $D^1_2$ and  $D^1_{2,3}$ , having codimension 3,
and their complement having codimension larger than or equal to 4.

Figures \ref{fig:d12sur} and \ref {fig:d23sur} illustrate the unfolding of
 umbilic  points in a generic one parameter family of immersions.
\end{proof}

\vskip .5cm

 \begin{remark} Although pictures equivalent to Fig. \ref{fig:d23sur} --
 $D^{1}_{2,3}$---
 appear in previous  printed work, no proofs seem to have been
 provided  for them before. See Gullstrand \cite{Gu} and  Porteous
 \cite{P3}.
 \end{remark}


\begin{remark}
There is a close relationship between one-parameter families of
immersed surfaces in $\mathbb R^3$ and immersions of 3-manifolds
in $\mathbb R^4$.

Assume that an  umbilic point $p$ at the value $0$ of the
parameter $\lambda$  in a family $\alpha_\lambda$ is non flat
--$k\ne 0$--.   By means of an inversion, which does not modify
the principal configuration, this can be achieved.

The principal configuration of the family around $(p,0)$ is the
same as that of the immersion $(u,v,\lambda)\to
(\alpha_\lambda(u,v), \lambda)$.

The generic  bifurcation patterns $D^1_2$ and $D^1_{2,3}$ studied
here correspond to the generic partially umbilic lines for
immersions of 3-manifolds in  $\mathbb R^4$. See Garcia
\cite{ga1}, \cite{ga2}.
\end{remark}

\section{Global Implications of Umbilic Bifurcations} \label{sec:3}

The local bifurcations at umbilic points discussed above have global implications on
the principal configurations of the immersions.  In order to formulate
 precisely the pertinent
 results it is
necessary to review some terminology established in previous papers.
See \cite {G-S1, G-S2, G-S3}

\subsection {Principal Cycles }  \label{3.1}

 A compact line  $ c$   of  ${\mathcal L}_{1\alpha}$ (resp. ${\mathcal L}_{2\alpha}$)
  is called minimal  (resp. maximal)  principal cycle of  $\alpha$.

Call $\pi = \pi_c$  the Poincar\'e first return map (holonomy) defined by the
lines of the foliation to which  $c$  belongs, defined on a segment of a line
of the orthogonal foliation through  a point  $o$  in  $c$.
A cycle is called  hyperbolic  if  $\pi^\prime(0)\ne 1$.  It has been proved
in \cite{G-S1} that
 $c$  is hyperbolic if and only if one of the following, equivalent
conditions hold

        a)  $\int_c dk_{2,\alpha}/(k_{2,\alpha}-k_{1,\alpha}) = \int_c dk_{1,\alpha}/(k_{2,\alpha}-k_{1,\alpha}) \ne 0$

       b)  $\int_c \frac {d{\mathcal H}_\alpha}{\sqrt{{\mathcal H}_{\alpha}^2-{\mathcal K}_{\alpha}}}\ne 0$.

The simplest bifurcation of non-hyperbolic principal cycles has been
studied in \cite{G-S3} and \cite{Ga-S2}.

\subsection{Umbilic Connections and Loops.\\}\label{3.2}

 A principal line  $\gamma$  which is an {\it umbilic  separatrix} of two
 different umbilic
points  $p$, $q$  of  $\alpha$ or twice a separatrix of the same
umbilic point $ p$  of $\alpha$ is called an {\it umbilic
separatrix connection} of  $\alpha$;  in the second case $\gamma$
is also called an {\it umbilic separatrix loop.} The simplest
bifurcations of umbilic connections as well as the consequent
appearance of principal cycles have been studied in \cite{G-S5}.

There are two bifurcation patterns producing  principal cycles which
will be studied in this work.  They are associated with the bifurcations
of  $D^{1}_{2}$  and  $D^{1}_{2,3}$  umbilic points, when their
separatrices form loops, self connecting these points. They are defined as
follows.

\begin{figure}[htbp]
 \begin{center}
 \includegraphics[angle=0, width=13cm]{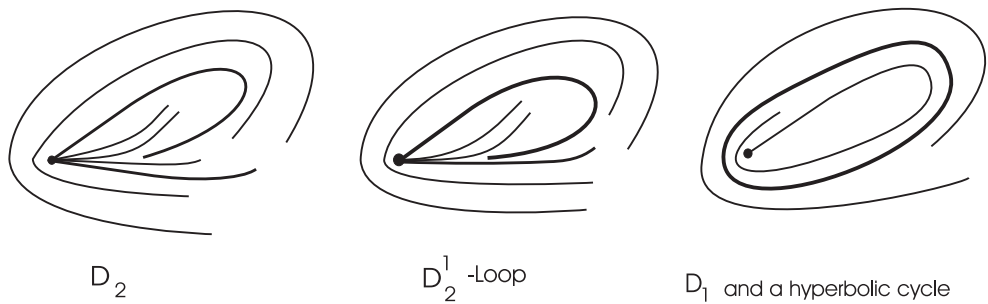}
 \caption{ $D^{1}_{2}$ -   loop bifurcation \label{fig:d12loop}}
   \end{center}
 \end{figure}

A $D^{1}_{2}$ - {\it interior loop} consists on a point of type $ D^{1}_{2}$ and its
isolated separatrix, which is assumed to be  contained in the interior of
the parabolic sector. See  Fig. \ref{fig:d12loop}
, where such loop together with the
bifurcating principal cycle are illustrated. Here, a hyperbolic principal
cycle bifurcates from the loop when the  $D^{1}_{2}$ point bifurcates into  $D_1$ .

\vskip .3cm

\begin{figure}[htbp]
 \begin{center}
 \includegraphics[angle=0, width=13cm]{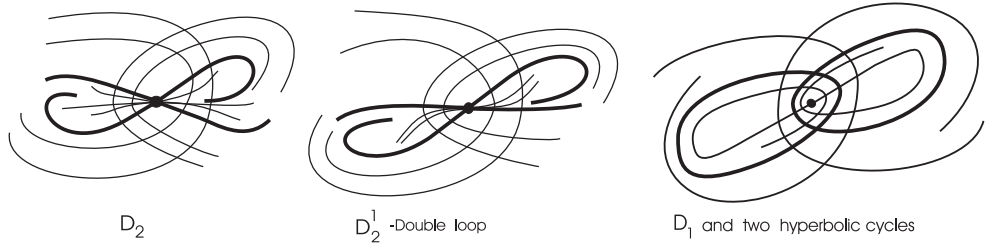}
 \caption{ $D^{1}_{2}$ -  double  loop bifurcation \label{fig:d12loopd}}
   \end{center}
 \end{figure}

    If both principal foliations have   $D^{1}_{2}$  - interior loops (at the same
 $D^{1}_{2}$   point), after bifurcation there appear two hyperbolic cycles, one for
each foliation. This case will be called {\it double} $D^{1}_{2}$
 - {\it interior loop}. In
Fig. \ref{fig:d12loopd},
Fig. \ref{fig:d12loop}
has been  modified and completed accordingly so as to
represent both maximal and minimal foliations, each  with its respective
 $D^{1}_{2}$ -interior loops (left) and bifurcating  hyperbolic principal
  cycles (right).

\begin{proposition}\label{pr:d12lup}
A unique hyperbolic principal cycle bifurcates from a $D^{1}_{2}$   loop.
\end{proposition}
\begin {proof}
The proof follows from the same analysis leading  to the uniqueness and hyperbolicity
of the periodic orbit bifurcating from  the singular cycle consisting of separatrix
 connecting a saddle-node   and a saddle through a  separatrix  located at
 projective line
   and a finite saddle separatrix, interior to the parabolic sector of the
   saddle-node. See \cite{PS}.

        \begin{figure}[htbp]
         \begin{center}
         \includegraphics[angle=0, width=5cm]{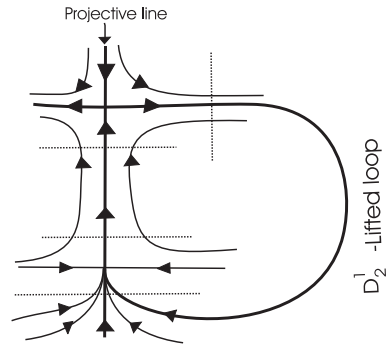}
         \caption{ $D^{1}_{2}$ - lifted loop. \label{fig:d12laco}}
           \end{center}
         \end{figure}

   This  consists in the decomposition of the
   return map into a singular transition (across the saddle-node point),
   whose contracting behavior is very small, even when compared with the
    transition along a saddle hyperbolic cycle,  and two regular transitions
    along the separatrix at infinity and a regular one (along the finite
    saddle separatrix), whose expansion behavior is bounded. See
    Fig.
    \ref{fig:d12laco}.
\end{proof}

\begin{figure}[htbp]
 \begin{center}
 \includegraphics[angle=0, width=12cm]{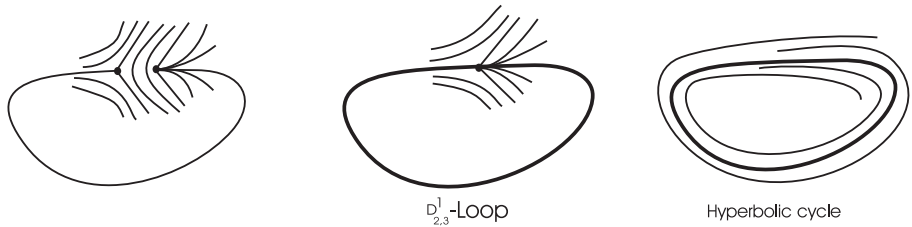}
 \caption{ $D^{1}_{2,3}$ - loop bifurcation. \label{fig:d23loop}}
   \end{center}
 \end{figure}

    A $D^{1}_{2,3}$ - {\it interior loop} consists on a point of type $D^{1}_{2,3}$
     and its
hyperbolic separatrix, which is assumed to be  contained in the interior of
the parabolic sector.
See  Fig. \ref{fig:d23loop}, where such loop together with the
bifurcating principal cycle are illustrated. Here, a hyperbolic principal
cycle bifurcates from the loop when the umbilic points are annihilated.

Notice that the $D^{1}_{2}$ and  $D^{1}_{2,3}$   interior loops described above, being
separatrices at only one and,   are neither umbilic
separatrix connections or loops in the sense of Gutierrez and  Sotomayor \cite{G-S5}.

 \begin{proposition}
A unique hyperbolic principal cycle bifurcates from a $D^{1}_{2,3}$   loop.
\end{proposition}
\begin {proof}
The proof follows from the same analysis leading to the uniqueness and hyperbolicity of
 the saddle-node loop bifurcation of vector fields on 2-manifolds \cite{S}. This
  consists in the decomposition of the return map into a singular transition
  (across the saddle-node point), whose contracting behavior is very small and
  a regular one (along the loop), whose expansion behavior is bounded.
\end{proof}

\section{ Relative Principal Structural Stability and Globalization.}\label{sec:4}

   Let $ \mathcal C$
be a subset of  ${\mathcal I}^{r,s}.$
 An element  $\alpha \in {\mathcal C}\subset {\mathcal I}^{r,s}$ is said $C^s-$
structurally stable relative to   $\mathcal C$,  if there is a neighborhood  $\mathcal V$  of  $\alpha$
in ${\mathcal I}^{r,s}$  such that for every    $\beta\in {\mathcal C}\cap{\mathcal  V}$
there is a homeomorphism  $ h =
h_\beta$  of $\mathbb M ^2$ which maps  ${\mathcal U}_\alpha$ onto  ${\mathcal U}_\beta$
 and maps lines of
${\mathcal F}_{i,\alpha}   $ onto    those of  ${\mathcal F}_{i,\beta}, \; i=1,2.   $

 When  $\mathcal C$ is the whole  ${\mathcal I}^{r,s}$,  $\alpha$ is
called simply $C^s$-structurally stable.

Call  ${\mathcal S}^r(j)$,  $j = a,b,c,d, $ the set
of  $\alpha \in  {\mathcal I}^{r}$, $  r \geq 4$  such  that,
respectively,
\begin{enumerate}
   \item[a)]  all the umbilic points of  $\alpha$ are  Darbouxian,

   \item[ b)] all the principal cycles of  $\alpha$ are hyperbolic,

    \item[c)] $\alpha$ has no umbilic separatrix connections,

   \item[ d)] the limit set of every principal line of  $\alpha$ is the union of
    umbilic points, principal cycles and umbilic connections.
\end{enumerate}
    The basic stability and genericity result that follows provides a
synthesis of these conditions.

\begin{theorem}    \label{th:gs}
    (\cite{G-S1,G-S2,G-S3}). Let $ r \geq 4.$
  The set  ${\mathcal S}^r= \cap_{\{j=a,\cdots,d\}} {\mathcal S}^r(j)$, 
  is open in   ${\mathcal I}^{r,3}$  and dense in
    ${\mathcal I}^{r,2}$ .   Every
     $\alpha \in  {\mathcal S}^r$    is $C^3$-structurally stable.
 \end{theorem}

The structural stability for lines of curvature has been extended
recently  to other equations and foliations of classical geometry.
See
   \cite{ggs}, \cite{arx}, \cite {axial}, \cite{amean},
  \cite{tou}, \cite{sbm},

    Recall that a {\it  Banach submanifold } of class $C^k$ and codimension one
of the Banach manifold  ${\mathcal I}^{r,r} =  {\mathcal I}^{r,r} ({\mathbb M}^2,{\mathbb R}^3)$
 is a subset $\mathcal B$ locally implicitly defined by
the set of zeroes of a real valued $C^k$  function with non
vanishing derivative, see \cite{S} and \cite{la}.

Now it is possible to state the main global results of this paper.

Call  ${\mathcal S}^{r}(a_1)$  the set of  $\alpha \in
 {\mathcal I}^{r }\backslash {\mathcal S}^{r}(a)$
  with a non-Darbouxian
umbilic point at which the transversality condition 
$T$  holds.

\begin{theorem}   \label{th:d12gl}
      Let $r \geq 5.$ Then the  following holds.
      \begin{enumerate}
      \item[i)]
 The set ${\mathcal S}_1^{r}(a_1)$   of immersions
 $\alpha \in \cap_{\{j=b,c,d\}}{\mathcal S}^{r}(j)$, 
  \-such
that   all their umbilic points are  Darbouxian
 except one which is of type $D^{1}_{2}$,
 is a Banach submanifold of codimension  1 and of class
$C^{r-3}$  of $ {\mathcal I}^{r,r}$.
 \item[ii)] The set
${\mathcal S}_1^{r}(a_1)$  is open in  $ {\mathcal I}^{r,4}=
   {\mathcal I}^{r}\backslash  {\mathcal S}^{r} $
  endowed with the $C^3$-topology, and is dense in
    ${\mathcal I}_1^{r,2}(a_1) = {\mathcal I}_1^{r}(a_1)$
      endowed with the $C^2$-topology.
      \item[iii)]
Every  $\alpha \in  {\mathcal S}_1^{r}(a_1)$  is
$C^4$-structurally stable relative to ${\mathcal S}^{r}(a_1)$.
\end{enumerate}
 \end{theorem}

\begin{proof} Outline.
Through the Lie-Cartan lifting, the proof is localized and reduced
to the case of vector fields where the Banach structure
established in the work of Sotomayor \cite {S} applies. In this
case the center manifold of the saddle-node is tangent to the
projective line and the characterization of the umbilic point
$D^1_2$ is in terms of the 3-jet of the immersion. The methods of
canonical regions and their continuous dependence of the
immersions used in \cite{G-S1} apply with minor modifications to
the present case to achieve the relative openness and construction
of homomorphism preserving principal configurations between
immersions inside ${\mathcal S}_1^{r}(a_1)$ which are close to
each other. The approximation leading to the $C^{r,2}$ density is
similar to that developed by Gutierrez and Sotomayor \cite{G-S2,
G-S3} to prove the $C^{r,2}$ density of ${\mathcal S}^r.$
\end{proof}

Call  ${\mathcal I}_1^{r}(a_2)$    the set of
 $\alpha \in  {\mathcal I}^{r }\backslash {\mathcal S}^{r}(a)$
  with an umbilic point at
which the transversality condition $T$  does not hold.

\newpage
    \begin{theorem} \label{th:d23gl}
Let $r \geq 5$. Then the following holds
\begin{enumerate}
\item[i)]
 The set ${\mathcal S}_1^{r}(a_2)$   of immersions
  $\alpha \in \cap_{\{j=b,c,d\}}{\mathcal S}^{r}(j)$, 
  such
that   all their umbilic points are  Darbouxian  except
 one which is of type $D^{1}_{2,3}$,
 is a Banach submanifold of codimension  1 and of class
$C^{r-3}$  of $ {\mathcal I}^{r,r}$.

\item[ii)] The set ${\mathcal S}_1^{r}(a_2)$  is open in $
{\mathcal I}^{r,5}= {\mathcal I}^{r}\backslash  {\mathcal S}^{r} $
  endowed with the $C^3$-topology, and is dense in
${\mathcal I}_1^{r,2}(a_2)= {\mathcal I}_1^{r}(a_2)$
    endowed with the $C^2$-topology.

\item[iii)] Every  $\alpha \in  {\mathcal S}_1^{r}(a_2)$ is
$C^5$-structurally stable relative to ${\mathcal S}^{r}(a_2)$.
\end{enumerate}
\end{theorem}

\begin{proof} Outline.
Through the Lie-Cartan lifting, the proof is localized and reduced
to the case of vector fields where the Banach structure
established in the work of Sotomayor \cite{S} applies. In this
case the center manifold of the saddle-node is transversal  to the
projective line and the characterization of the umbilic point
$D^1_{2,3}$ is in terms of the 4-jet of the immersion. The methods
of canonical regions and their continuous dependence of the
immersions used in \cite{G-S1} apply with minor modifications to
the present case to achieve the relative openness and construction
of homomorphism preserving principal configurations between
immersions inside ${\mathcal S}_1^{r}(a_2)$ which are close to
each other. The approximation leading to the $C^{r,2}$ density is
similar to that developed by Gutierrez and Sotomayor \cite{G-S2,
G-S3} to prove the $C^{r,2}$ density of ${\mathcal S}^r.$
\end{proof}

\section {Concluding Remark}\label{cr}

In this work  we have established that
     Figures \ref{fig:d12bifu} or \ref{fig:d12sur} and \ref{fig:bif23} or
     \ref{fig:d23sur}, completed when pertinent  with the global
implications in Figs.
     \ref{fig:d12loop}, \ref{fig:d12loopd} and \ref{fig:d23loop},
     represent the
topological changes on  the principal configurations
 of families of  immersions   depending on a parameter. These
bifurcations are  exhibited when the families
cross  transversally the codimension one submanifolds ${\mathcal
S}_1^{r}(a_1)$
and ${\mathcal S}_1^{r}(a_2)$. This contribution, together with those
presented in  \cite{G-S4} and \cite{G-S5},  completes the
analysis of
codimension one  bifurcations for principal configurations formulated in
 \cite{G-S7}.

 \newpage

\author{\noindent Carlos Guti\'errez\\Instituto de Ci\^encias  Matem\'{a}ticas
e Computa\c c\~ao,\\Universidade de S\~{a}o Paulo,  Caixa Postal 668,
\\CEP 13560-970, S\~{a}o Carlos, S.P., Brazil \\
\\
Jorge Sotomayor\\Instituto de Matem\'{a}tica e Estat\'{\i}stica,\\
Universidade de S\~{a}o Paulo, \\Rua do Mat\~{a}o 1010, Cidade Universit\'{a}ria,
 \\CEP 05508-090, S\~{a}o Paulo, S.P., Brazil \\
\\ Ronaldo Garcia\\Instituto de Matem\'{a}tica e Estat\'{\i}stica,\\
Universidade Federal de Goi\'as,\\CEP 74001-970, Caixa Postal 131,\\
Goi\^ania, GO, Brazil}

\end{document}